\documentclass[a4paper,12pt]{article}
\usepackage{amsfonts,amsmath,amssymb}

\makeatletter
\newbox\bk@bxb
\newbox\bk@bxa
\newif\if@bkcont
\newcount\bk@lcnt

\def\breakboxskip{2pt}
\def\breakboxparindent{1.8em}

\def\breakbox{\vskip\breakboxskip\relax
\setbox\bk@bxb\vbox\bgroup
\advance\linewidth -2\fboxrule
\hsize\linewidth\@parboxrestore
\parindent\breakboxparindent\relax}

\def\bk@split{%
\@tempdimb\ht\bk@bxb 
\advance\@tempdimb\dp\bk@bxb
\setbox\bk@bxa\vsplit\bk@bxb to\z@ 
\setbox\bk@bxa\vbox{\unvbox\bk@bxa}
\setbox\@tempboxa\vbox{\copy\bk@bxa\copy\bk@bxb}
\advance\@tempdimb-\ht\@tempboxa
\advance\@tempdimb-\dp\@tempboxa}

\def\bk@addfsepht{%
\setbox\bk@bxa\vbox{\vskip\fboxsep\box\bk@bxa}}

\def\bk@addskipht{%
\setbox\bk@bxa\vbox{\vskip\@tempdimb\box\bk@bxa}}

\def\bk@addfsepdp{%
\@tempdima\dp\bk@bxa
\advance\@tempdima\fboxsep
\dp\bk@bxa\@tempdima}

\def\bk@addskipdp{%
\@tempdima\dp\bk@bxa
\advance\@tempdima\@tempdimb
\dp\bk@bxa\@tempdima}

\def\bk@line{%
\hbox to \linewidth{%
\hskip-2\fboxsep\vrule \@width\fboxrule\hskip.5\fboxsep\vrule \@width\fboxrule\hskip1.5\fboxsep
\box\bk@bxa\hfil
}}%

\def\endbreakbox{\egroup
\ifhmode\par\fi{\noindent\bk@lcnt\@ne
\@bkconttrue\baselineskip\z@\lineskiplimit\z@
\lineskip\z@\vfuzz\maxdimen
\bk@split\bk@addfsepht\bk@addskipdp
\ifvoid\bk@bxb 
\def\bk@fstln{\bk@addfsepdp
\hskip-\parindent\vbox{\llap{\raisebox{-2ex}{\rule{1.5\fboxsep}{\fboxrule}\hskip.5\fboxsep}}\bk@line\llap{\rule{1.5\fboxsep}{\fboxrule}\hskip.5\fboxsep}}}

\else 
\def\bk@fstln{\vbox{\llap{\raisebox{-2ex}{\rule{1.5\fboxsep}{\fboxrule}\hskip.5\fboxsep}}\bk@line}\hfil%
\advance\bk@lcnt\@ne
\loop
\bk@split\bk@addskipdp\leavevmode
\ifvoid\bk@bxb 
\@bkcontfalse\bk@addfsepdp
\vtop{\bk@line\noindent\hskip-2\fboxsep{\rule{1.5\fboxsep}{\fboxrule}}}%

\else 
\bk@line
\fi
\hfil\advance\bk@lcnt\@ne
\if@bkcont\repeat}%
\fi
\leavevmode\bk@fstln\par}\vskip\breakboxskip\relax}

\newcommand{\scalv}[2]{\langle#1,#2\rangle}
\def\ret{\hspace{-.2ex}}
\def\smp{\smallskip\par}
\def\un{{\bf 1}}
\def\zero{\{0\}}
\def\pf{\noindent{\bf Proof~:}\ }
\def\findemo{~\leaders\hbox to 1em{\hss\  \hss}\hfill~\raisebox{.5ex}{\framebox[1ex]{}}\smp}
\def\spn{\bigskip\par\noindent}
\def\mpn{\medskip\par\noindent}
\def\smpn{\smallskip\par\noindent}
\def\normal{\mathop{\trianglelefteq}}

\def\mpoint{\;\;.}
\def\mvirg{\;\;,}

\def\Res{{\rm Res}}
\def\Ind{{\rm Ind}}
\def\Inf{{\rm Inf}}

\def\Hom{{\rm Hom}}

\def\Inf{{\rm Inf}}

\def\Ker{{\rm Ker}}

\def\Irr{{\rm Irr}}

\def\dsp{\displaystyle}
\def\Z{\mathbb{Z}}

\def\Q{\mathbb{Q}}

\newcommand{\romain}[1]{\uppercase\expandafter{\romannumeral #1}}
\newcommand{\gMod}[1]{#1{\hbox{-}\mathsf{Mod}}}

\newcommand{\sur}[1]{\,\overline{\! #1}}
\newcommand{\sumb}[2]{\mathop{\sum}_{{\scriptstyle #1}\atop {\scriptstyle #2}}}

\newcommand{\sumc}[3]{\sum_{{\scriptstyle #1}\atop {{\scriptstyle #2}\atop {\scriptstyle #3}}}}

\newenvironment{enonce}[1]{\pagebreak[2]\refstepcounter{subsection}\refstepcounter{prop}\smpn{{\bf \thesection.\arabic{prop}.\ \ #1~:}}\begin{it} }{\end{it}\smp}
\newenvironment{enonce*}[1]{\pagebreak[2]\smpn{#1~:}\begin{it} }{\end{it}\smp}
\newcommand{\result}[1]{\begin{enonce}{#1}}
\def\fresult{\end{enonce}}
\newcommand{\npar}{\smallskip\par\noindent\pagebreak[2]\refstepcounter{subsection}\refstepcounter{prop}{\bf \thesection.\arabic{prop}.\ \ }}

\newenvironment{mth}[1]{\pagebreak[3]\begin{breakbox}\begin{enonce}{#1}}{\end{enonce}\end{breakbox}}
\newenvironment{mth*}[1]{\pagebreak[3]\begin{breakbox}\begin{enonce*}{#1}}{\end{enonce*}\end{breakbox}}
\newenvironment{rem}[1]{\refstepcounter{subsection}\refstepcounter{prop} \mpn{{\bf \thesection.\arabic{prop}.}\ \ \bf#1\ :}}{\smp}

\def\dom{\backslash}
\renewenvironment{equation}{\refstepcounter{subsection}\refstepcounter{prop}$$}{\leqno{\bf (\theprop)}$$}

\newcommand{\sscalv}[2]{\langle\!\langle#1,#2\rangle\!\rangle}
\def\comack{\mathsf{coMack}}
\begin{document}
\centerline{\Large\bf On the Cartan matrix of Mackey algebras}\vspace{.5cm}
\centerline{Serge Bouc}\vspace{1cm}
{\footnotesize{\bf Abstract~:} Let $k$ be a field of characteristic $p>0$, and $G$ be a finite group. The first result of this paper is an explicit formula for the determinant of the Cartan matrix of the Mackey algebra $\mu_k(G)$ of $G$ over $k$. The second one is a formula for the rank of the  Cartan matrix of the cohomological Mackey algebra $co\mu_k(G)$ of $G$ over $k$, and a characterization of the groups $G$ for which this matrix is non singular. The third result is a generalization of this rank formula and characterization to blocks of $co\mu_k(G)$~: in particular, if $b$ is a block of $kG$, the Cartan matrix of the corresponding block $co\mu_k(b)$ of $co\mu_k(G)$ is non singular if and only if $b$ is nilpotent with cyclic defect groups.
\medskip\par
{\bf AMS Subject Classification~:} 18G05, 20C20, 20J06.\par
{\bf Keywords~:} Cartan matrix, cohomological, Mackey functor, block.
}
\section{Introduction}
\npar The theory of Mackey functors for a finite group $G$ over a commutative ring~$R$ is by many aspects very similar to the theory of $RG$-modules. It has been shown by J.~Th\'evenaz and P.~Webb (\cite{thevwebb}), among many fundamental other results, that the category of Mackey functors for $G$ over $R$ is equivalent to the category of modules over {\em the Mackey algebra} $\mu_R(G)$. This algebra shares many properties with the group algebra $RG$~: it is free as an $R$-module, and its $R$-rank does not depend on $R$~; if $K$ is a field of characteristic 0 or coprime to the order of $G$, the algebra $\mu_K(G)$ is semisimple~; when $(K,\mathcal{O},k)$ is a $p$-modular system, there is a decomposition theory from Mackey functors for $G$ over $K$ to Mackey functors for $G$ over $k$~; the Cartan matrix of $\mu_k(G)$ is symmetric and non singular. This list of common properties between the Mackey algebra and the group algebra is far from exhaustive\ldots\par
\npar However some well known results for group algebras are no longer true for Mackey algebras. It was observed on small examples in particular by M.~Nicollerat in her thesis (\cite{nicollerat} Chapitre 5) that the determinant of the Cartan matrix of $\mu_k(G)$, when $k$ is a field of characteristic $p$, is generally not a power of~$p$, even when $G$ itself is a $p$-group. Instead, some rather strange prime factors appear in this determinant. One of the motivations of this paper is to give an explanation for these strange factors, by stating an explicit formula for the determinant of the Cartan matrix of~$\mu_k(G)$ (Theorem~\ref{theorem Mackey}).\par
\npar The other motivation is the similar problem for the {\em cohomological} Mackey algebra $co\mu_k(G)$, which was defined by Th\'evenaz and Webb as a specific quotient of $\mu_k(G)$, with the property that the modules over $co\mu_k(G)$ are exactly the {\em cohomological Mackey functors}~: the first major difference is that in general, the Cartan matrix of this algebra is singular. This raises the question of characterizing those finite groups $G$ for which the Cartan matrix of $co\mu_k(G)$ is non singular, and the answer is the second main result of this paper (Theorem \ref{theorem coMackey})~: these groups are exactly the $p$-nilpotent groups with cyclic Sylow $p$-subgroups. The other possibly interesting result is that for an arbitrary finite group $G$, the rank of the Cartan matrix of $co\mu_k(G)$ is equal to the number of conjugacy classes of pairs $(R,s)$, where $R$ is a cyclic $p$-subgroup of $G$, and $s$ is a $p'$-element of the centralizer of $R$ in $G$. \par
The third result of this paper is a natural generalization of Theorem~\ref{theorem coMackey} to {\em blocks}, suggested by the one to one correspondence $b\mapsto co\mu_k(b)$ between blocks of $kG$ and blocks of $co\mu_k(G)$~: Theorem~\ref{theorem coMackey blocks} gives a formula for the rank of the Cartan matrix of $co\mu_k(b)$, in terms of $b$-Brauer pairs, and shows that this matrix is non singular if and only if the block $b$ is nilpotent with cyclic defect groups.
\npar The paper is organized according to these results~: Section~\ref{Mackey} is devoted to the case of $\mu_k(G)$, starting by recalling some standard notation, definitions and properties, and Section~\ref{coMackey} deals with $co\mu_k(G)$. Finally, Section~\ref{coMackey blocks} is devoted to the case of blocks of cohomological Mackey functors.\spn
{\bf Acknowledgment~:} I wish to thank Jacques Th\'evenaz for very fruitful discussions and friendly collaboration on all these questions also, at the EPFL in June~2009. 
\section{The Mackey algebra}\label{Mackey}
\npar Throughout the paper when $G$ is a finite group and $p$ is a prime number, the set of $p'$-elements of $G$ is denoted by $G_{p'}$, and the symbol $[G_{p'}]$ denotes a set of representatives of $G$-conjugacy classes in $G_{p'}$. The set of $p$-subgroups of~$G$ is denoted by $\mathcal{S}_p(G)$, and $[\mathcal{S}_p(G)]$ denotes similarly a set of representatives of $G$-conjugacy classes in $\mathcal{S}_p(G)$. 
\npar From now on $G$ will be a fixed finite group, and $p$ a fixed prime number. Let $(K,\mathcal{O},k)$ be a $p$-modular system~: thus $\mathcal{O}$ is a discrete valuation ring with residue field $k$ of characteristic $p$, and field of fractions $K$ of characteristic~0. Assume that $K$ and $k$ are splitting fields for all the groups $N_G(Q)/Q$, for $Q\in\mathcal{S}_p(G)$ (e.g. assume that $K$ contains the $e$-th roots of unity, where $e$ is the exponent of~$G$).\par
\npar When $R$ is a commutative ring with identity element, the Mackey algebra $\mu_R(G)$ of $G$ over $R$ has been defined by Th\'evenaz and Webb (\cite{thevwebb} Section~3). It is an associative algebra with the property that the category $\gMod{\mu_R(G)}$ of left $\mu_R(G)$-modules is equivalent to the category $\mathsf{Mack}_R(G)$ of Mackey functors for $G$ over $R$.
\npar Th\'evenaz and Webb have also shown (\cite{thevwebb} Theorem~10.1) that there is an equivalence of abelian categories
$$\mathsf{Mack}_k(G)\cong\prod_H\mathsf{Mack}_k\big(N_G(H)/H,\un\big)\mvirg$$
where $H$ runs through a set of representatives of conjugacy classes of $p$-perfect subgroups of $G$, and $\mathsf{Mack}_k\ret\big(\ret N_G(H)/H,\un\ret\big)$ denotes the subcategory of $\mathsf{Mack}_k\ret\big(\ret N_G(H)/H\ret\big)$ consisting of Mackey functors which are {\em projective relative} to $p$-subgroups. The category $\mathsf{Mack}_k(G,\un)$ is equivalent to $\gMod{\mu_k(G,\un)}$, where $\mu_k(G,\un)$ is a direct summand of $\mu_k(G)$ of the form $\mu_k(G)f$, for a specific central idempotent $f$ of $\mu_k(G)$.
\npar (\cite{thevwebb} Theorem 12.7 and Corollary 12.8) The correspondence $M\mapsto M(\un)$, that is evaluation at the trivial subgroup of $G$, induces a one to one correspondence between the set of isomorphism classes of indecomposable projective $\mu_k(G,\un)$-modules and the set of isomorphism classes of indecomposable $p$-permutation $kG$-modules (also called trivial source $kG$-modules).
\npar Let $pp_k(G)$ denote {\em the Green ring of $p$-permutation $kG$-modules} : as a group, it is the Grothendieck group of the category of finitely generated $p$-permutation $kG$-modules, for relations given by direct sum decompositions. The product on $pp_k(G)$ is induced by the tensor product of $kG$-modules over~$k$. \par
If $W$ is a finitely generated $p$-permutation $kG$-module, its dual $W^*=\Hom_k(W,k)$ is also a $p$-permutation $kG$-module, and this duality extends to a ring automorphism $W\mapsto W^*$ of $pp_k(G)$. \par
When $Q$ is a $p$-subgroup of $G$, the {\em Brauer quotient} $W[Q]$ is a $p$-permutation $k\sur{N}_G(Q)$-module, where $\sur{N}_G(Q)=N_G(Q)/Q$ (see \cite{brouepperm}). This construction commutes with duality and tensor product of $p$-permutation modules~: if $V$ and $W$ are (finitely generated) $p$-permutation $kG$-modules, there are isomorphisms of $kN_G(Q)/Q$-modules 
\begin{equation}\label{iso}
W[Q]^*\cong W^*[Q],\;\;\;V[Q]\otimes W[Q]\cong (V\otimes W)[Q]\mpoint
\end{equation}
In particular, the Brauer quotient induces a ring homomorphism from $pp_k(G)$ to $pp_k\big(\sur{N}_G(Q)\big)$, still denoted by $W\mapsto W[Q]$ (see e.g. \cite{both5} Proposition~2.11).\par
There is a $\Z$-valued bilinear form $\scalv{~\,}{~}$ on $pp_k(G)$ defined for $p$-permutation $kG$-modules $V$ and $W$ by
$$\scalv{V}{W}_G=\dim_k\Hom_{kG}(V,W)\mpoint$$
It is worth noticing that
\begin{equation}\label{scalprod}
\scalv{V}{W}_G=\dim_k\Hom_{kG}(k,V^*\otimes_kW)\mpoint
\end{equation}
The bilinear form $\scalv{~\,}{~}_G$ extends to a $K$-valued bilinear form on $K\otimes_\Z pp_k(G)$, still denoted by $\scalv{~\,}{~}_G$.
\npar \label{sscalv}Let $\sscalv{~\,}{~}_G$ denote the bilinear form on $pp_k(G)$ defined for $p$-permutation $kG$-modules $V$ and $W$ by
\begin{eqnarray*}
\sscalv{V}{W}_G&=&\sum_{Q\in[\mathcal{S}_p(G)]}\scalv{V[Q]}{W[Q]}_{\sur{N}_G(Q)}\\
&=&\sum_{Q\in[\mathcal{S}_p(G)]}\dim_k\Hom_{k\sur{N}_G(Q)}(V[Q],W[Q])\mpoint
\end{eqnarray*} 
It follows from \ref{iso} that
\begin{equation}\label{scalaire}
\sscalv{V}{W}_G=\sscalv{k}{V^*\otimes_kW}_G\mpoint
\end{equation}
It was shown in \cite{resomackey} Proposition~5.11 that if $L$ and $M$ are projective Mackey functors in $\mathsf{Mack}_k(G,\un)$, then
$$\dim_k\Hom_{\mathsf{Mack}_k(G,\un)}(L,M)=\sscalv{L(\un)}{M(\un)}_G\mpoint$$
When $L$ and $M$ are indecomposable, this is equal to the coefficient $c_{L,M}$ of the Cartan matrix of the algebra $\mu_k(G,\un)$.
\npar Let $\mathcal{Q}_{G,p}$ denote the set of pairs $(R,s)$ consisting of a $p$-subgroup $R$ of~$G$ and a $p'$-element $s$ of $\sur{N}_G(R)$. The group $G$ acts by conjugation on $\mathcal{Q}_{G,p}$. Let $[\mathcal{Q}_{G,p}]$ denote a set of representatives of $G$-orbits on $\mathcal{Q}_{G,p}$, and let $N_G(R,s)$ denote the stabilizer of $(R,s)\in\mathcal{Q}_{G,p}$ in $G$~: it is the set of elements $g\in N_G(R)$ such that the image $\sur{g}$ of $g$ in $\sur{N}_G(R)$ centralizes $s$. In other words, there is an exact sequence of groups 
\begin{equation}\label{suite exacte}
\un\to R\to N_G(R,s)\to C_{\sur{N}_G(R)}({s})\to \un\mpoint
\end{equation}
It was shown in \cite{both5} that {\em the primitive idempotents} of the (commutative) ring $K\otimes_\Z pp_k(G)$ are indexed by the orbits of $G$ on $\mathcal{Q}_{G,p}$~: the idempotent $F_{R,{s}}^G$ associated to the orbit of the pair $(R,s)$ is equal to
\begin{equation}\label{idempotent}F_{R,{s}}^G=\frac{1}{|R||{s}||C_{\sur{N}_G(R)}({s})|}\sum_{\varphi,L}\tilde{\varphi}({s}^{-1})|L|\mu(L,{<}sR{>})\,\Ind_L^G\Res_L^{{<}sR{>}}k_\varphi\mvirg
\end{equation}
where 
\begin{itemize}
\item the group ${<}sR{>}$ is the inverse image in $N_G(R)$ of the subgroup ${<}s{>}$ of $\sur{N}_G(R)$ under the map $x\mapsto xR$.
\item the morphism $\varphi$ runs through group homomorphisms ${<}{s}{>}\to k^\times$, and $\tilde{\varphi}$ lifts $\varphi$ to $K$. The module $k_\varphi$ is the vector space $k$ on which ${<}sR{>}$ acts by ${<}sR{>}\to {<}s{>}\stackrel{\varphi}{\to}k^\times$ (thus $\tilde{\varphi}$ is the Brauer character of the module $k_\varphi$).
\item the group $L$ runs through the set of subgroups of ${<}sR{>}$ such that $LR={<}sR{>}$.
\end{itemize}
The idempotents $F_{R,{s}}^G$, for $(R,{s})\in[\mathcal{Q}_{G,p}]$, form a $K$-basis of $K\otimes_\Z pp_k(G)$. Any element $W$ of $K\otimes_\Z pp_k(G)$ can be expressed in this basis as
\begin{equation}\label{coeffs}W=\sum_{(R,{s})\in[\mathcal{Q}_{G,p}]}F_{R,{s}}^G\otimes_kW=\sum_{(R,{s})\in[\mathcal{Q}_{G,p}]}t_{R,{s}}^G(W)\,F_{R,{s}}^G\mvirg
\end{equation}
where $t_{R,{s}}^G$ is the extension to $K\otimes_\Z pp_k(G)$ of the {\em species} (i.e. the ring homomorphism, see \cite{benson-lnm} Lemma 2.2.1, page 26) $pp_k(G)\to K$ associated to the pair $(R,{s})$.
\npar\label{Brauer trace} Recall that the value of $t_{R,{s}}^G$ on a $p$-permutation $kG$-module $W$ is equal to the value at $s$ of the Brauer character of the Brauer quotient $W[R]$ (see \cite{both5} Notation~2.15 and Proposition 2.18). This will be also be called the {\em Brauer trace} of $s$ on $W[R]$, and denoted by ${\rm BrTr}(s\mid W[R])$.
\npar The ring automorphism $W\mapsto W^*$ extends by $K$-linearity to an automorphism of the ring $K\otimes_\Z pp_k(G)$. This automorphism preserves the set of primitive idempotents. The bilinear form $\sscalv{~\,}{~}_G$ also extends to a $K$-valued bilinear form on $K\otimes_\Z pp_k(G)$. Since the dual of $k_\varphi$ is isomorphic to $k_{\varphi^{-1}}$, and since $\widetilde{(\varphi^{-1})}({s}^{-1})=\tilde{\varphi}({s})$, it follows from \ref{idempotent} that $(F_{R,{s}}^G)^*=F_{R,{s}^{-1}}^G$. Now Equation~\ref{scalaire} shows that
$$\sscalv{F_{R',{s}'}^G}{F_{R,{s}}^G}_G=\left\{\begin{array}{cl}
\sscalv{k}{F_{R,{s}}^G}_G&\hbox{if}\;(R',{s}')=_G(R,{s}^{-1})\\
0&\hbox{otherwise}\end{array}\right.$$
(where $=_G$ denotes $G$-conjugacy). Thus, if $V$ and $W$ are any elements of $K\ret\otimes_\Z\nolinebreak pp_k(G)$, it follows from \ref{coeffs} that
$$\sscalv{V}{W}_G=\sum_{(R,{s})\in[\mathcal{Q}_{G,p}]}t_{R,{s}^{-1}}^G(V)t_{R,{s}}^G(W)\sscalv{k}{F_{R,{s}}^G}_G\mpoint$$
\npar {\em The indecomposable $p$-permutation $kG$-modules} are indexed by conjugacy classes of pairs $(P,E)$, where $P$ is a $p$-subgroup of $G$, and $E$ is an indecomposable projective $k\sur{N}_G(P)$-module (\cite{brouepperm} Theorem~3.2). Let $[\mathcal{P}_{G,p}]$ denote a set of representatives of such conjugacy classes in $G$. \par
The indecomposable module $M_{P,E}$ indexed by the pair $(P,E)$ has vertex $P$, and the projective module $E$ is isomorphic to the Brauer quotient $M_{P,E}[P]$. In particular, if $Q$ is a subgroup of $G$, then $M_{P,E}[Q]$ is non zero only if $Q\leq_GP$ (and in fact, if and only if $Q\leq_G P$). 
\npar The modules $M_{P,E}$, for $(P,E)\in [\mathcal{P}_{G,p}]$, form a basis of $pp_k(G)$, hence also a basis of $K\otimes_\Z pp_k(G)$. The Cartan matrix $\mathsf{C}\big(\mu_k(G,1)\big)$ of the algebra $\mu_k(G,\un)$ can be viewed as the square matrix indexed by $[\mathcal{P}_{G,p}]$, where the coefficient indexed by the elements $(P,E)$ and $(Q,F)$ of $[\mathcal{P}_{G,p}]$ is equal to
\begin{eqnarray*}
c_{(P,E),(Q,F)}&=&\sscalv{M_{P,E}}{M_{Q,F}}_G\\
&=&\sum_{(R,{s})\in[\mathcal{Q}_{G,p}]}t_{R,{s}^{-1}}^G(M_{P,E})t_{R,{s}}^G(M_{Q,F})\sscalv{k}{F_{R,{s}}^G}_G\mpoint
\end{eqnarray*}
\npar \label{a suivre}Let $\mathsf{T}$ and $\mathsf{T}'$ denote the matrices indexed by the product $[\mathcal{P}_{G,p}]\times\nolinebreak[4] [\mathcal{Q}_{G,p}]$, defined by
\begin{eqnarray*}
\mathsf{T}_{(P,E),(R,{s})}&=&t_{R,{s}}^G(M_{P,E})\mvirg\\
\mathsf{T}'_{(P,E),(R,{s})}&=&t_{R,{s}^{-1}}^G(M_{P,E})\mvirg
\end{eqnarray*}
respectively. Note that $\mathsf{T}$ and $\mathsf{T}'$ are square matrices. Let moreover $\mathsf{S}$ be the (square) diagonal matrix indexed by $[\mathcal{Q}_{G,p}]$, where the diagonal term indexed by $(R,{s})$ is equal to $\sscalv{F_{R,{s}}^G}{k}_G$. It follows that
$$\mathsf{C}\big(\mu_k(G,\un)\big)=\mathsf{T}'\cdot\mathsf{S}\cdot{^t\mathsf{T}}\mpoint$$
Thus
\begin{equation}\label{detproduit}
\det \mathsf{C}\big(\mu_k(G,\un)\big)=\det \mathsf{T}'\,\det\mathsf{S}\,\det {^t\mathsf{T}}\mvirg
\end{equation}
and moreover, since $\mathsf{S}$ is diagonal
\begin{equation}\label{detS}
\det\mathsf{S}=\prod_{(R,{s})\in[\mathcal{Q}_{G,p}]}\sscalv{k}{F_{R,{s}}^G}_G=\prod_{R\in[\mathcal{S}_p(G)]}\;\prod_{{s}\in[\sur{N}_G(R)_{p'}]}\sscalv{k}{F_{R,{s}}^G}_G
\end{equation}
\npar Now since $M_{P,E}[R]=\zero$ unless $R\leq_G P$, the matrices $\mathsf{T}$ and $\mathsf{T}'$ are actually block triangular, with diagonal blocks $\Delta(P)$ indexed by the conjugacy classes of $p$-subgroups $P$ of~$G$. Thus
$$\det \mathsf{T}=\prod_{P\in[\mathcal{S}_p(G)]}\det\Delta(P)\mpoint$$
The block matrix $\Delta(P)$ is a (square) matrix with rows indexed by the isomorphism classes of indecomposable projective $k\sur{N}_G(P)$-modules $E$, and columns indexed by the conjugacy classes of $p'$-elements ${s}$ of $\sur{N}_G(P)$. The coefficient $\Delta(P)_{E,{s}}$ is equal to
\begin{eqnarray*}
\Delta(P)_{E,{s}}&=&t_{P,{s}}^G(M_{P,E})\\
&=&{\rm Brauer\, trace}({s}\mid M_{P,E}[P])\\
&=&\Phi_E({s})\mvirg
\end{eqnarray*}
since $M_{P,E}[P]\cong E$, where $\Phi_E$ is the Brauer character of the projective $k\sur{N}_G(P)$-module $E$. \par
Similarly,
$$\det \mathsf{T}'=\prod_{P\in[\mathcal{S}_p(G)]}\det\Delta'(P)\mvirg$$
where the coefficients of the diagonal block $\Delta'(P)$ of $\mathsf{T}'$ are given by
$$\Delta'(P)_{E,{s}}=\Phi_E({s}^{-1})\mpoint$$

\npar Let $\Sigma(P)$ denote the (square) diagonal matrix indexed by the conjugacy classes of $p'$-elements of $\sur{N}_G(P)$, with diagonal coefficient $\Sigma(P)_{{s},{s}}$ equal to the inverse of the order of the centralizer $C_{\sur{N}_G(P)}({s})$ of ${s}$ in $\sur{N}_G(P)$, and let $U(P)=\Delta(P)\cdot\Sigma(P)\cdot {^t\Delta'(P)}$. The coefficient of $U(P)$ indexed by the indecomposable projective $k\sur{N}_G(P)$-modules $E$ and $F$ is equal to
\begin{eqnarray*}
U(P)_{E,F}&=&\sum_{{s}\in[\sur{N}_G(P)_{p'}]}\Phi_E({s})\frac{1}{|C_{\sur{N}_G(P)}({s})|}\Phi_F({s}^{-1})\\
&=&\frac{1}{|\sur{N}_G(P)|}\sum_{{s}\in \sur{N}_G(P)_{p'}}\Phi_E({s})\Phi_F({s}^{-1})\\
&=&\scalv{E}{F}_{\sur{N}_G(P)}\mpoint
\end{eqnarray*}
In other words, the matrix $U(P)$ is equal to the Cartan matrix of the group algebra $k\sur{N}_G(P)$. It is then well known (see \cite{brauer-nesbitt} \romain{3}.16 or \cite{feit}) that
$$\det U(P)=\prod_{{s}\in[\sur{N}_G(P)_{p'}]}|C_{\sur{N}_G(P)}({s})|_p\mvirg$$
and it follows that
\begin{eqnarray*}
\det \Delta(P)\det{^t\Delta'(P)}&=&\frac{\dsp\prod_{{s}\in[\sur{N}_G(P)_{p'}]}|C_{\sur{N}_G(P)}({s})|_p}{\det\Sigma(P)}\\
&=&\prod_{{s}\in[\sur{N}_G(P)_{p'}]}|C_{\sur{N}_G(P)}({s})|_p\prod_{{s}\in[\sur{N}_G(P)_{p'}]}|C_{\sur{N}_G(P)}({s})|\mpoint
\end{eqnarray*}
Together with \ref{detproduit} and \ref{detS}, this gives
\begin{equation}\label{detmu}
\det\mathsf{C}\big(\mu_k(G,\un)\big)=\prod_{R\in[\mathcal{S}_p(G)]}\prod_{{s}\in[\sur{N}_G(R)_{p'}]}\big(|C_{\sur{N}_G(R)}({s})|_p|C_{\sur{N}_G(R)}({s})|\sscalv{k}{F_{R,{s}}^G}_G\big)\mpoint
\end{equation}
\begin{mth}{Lemma} \label{Brauer F}Let $Q\in\mathcal{S}_p(G)$ and $(R,{s})\in\mathcal{Q}_{G,p}$. Then
$$F_{R,{s}}^G[Q]=\sumb{x\in N_G(Q)\dom G/N_G(R,s)}{Q\normal {{^x\!<}Rs{>}}}F_{^xR/Q,{^x{s}}}^{\sur{N}_G(Q)}\mvirg$$
(where $^x{s}$ denotes the image in $\sur{N}_G(Q,{^xR})/{^xR}$ of the $x$-conjugate of ${s}$).
\end{mth}
\pf The Brauer map $W\mapsto W[Q]$ is a ring homomorphism from $pp_k(G)$ to $pp_k\big(\sur{N}_G(Q)\big)$. Thus $F_{R,{s}}^G[Q]$ is an idempotent of $pp_k\big(\sur{N}_G(Q)\big)$, hence a sum of distinct primitive idempotents $F_{X,u}^{\sur{N}_G(Q)}$ associated to some pairs $(X,u)\in\mathcal{Q}_{\sur{N}_G(Q),p}$. In other words $X$ is a $p$-subgroup of $\sur{N}_G(Q)$, of the form $Y/Q$, for some $p$-subgroup~$Y$ of $G$ such that $Q\normal Y$, and $u$ is a $p'$-element of $\sur{N}_{\sur{N}_G(Q)}(X)\cong N_G(Q,Y)/Y$. \par
The idempotent $F_{X,u}^{\sur{N}_G(Q)}$ appears in $F_{R,{s}}^G[Q]$ if and only if the species $t_{X,u}^{\sur{N}_G(Q)}$ takes the value 1 when evaluated at $F_{R,{s}}^G[Q]$. But $t_{X,u}^{\sur{N}_G(Q)}(F_{R,{s}}^G[Q])$ is equal to the Brauer trace (see~\ref{Brauer trace}) of $u$ on $F_{R,{s}}^G[Q][Y/Q]=\Res_{N_G(Q,Y)/Y}^{N_G(Y)/Y}F_{R,{s}}^G[Y]$. In other words
$$t_{X,u}^{\sur{N}_G(Q)}(F_{R,{s}}^G[Q])=t_{Y,{u}}^G(F_{R,s}^G)\mvirg$$
and this is equal to 1 if and only if the pair $(Y,u)$ is $G$-conjugate to the pair $(R,s)$. The lemma follows.\findemo
\pagebreak[4]
\begin{mth}{Lemma} \label{scal F}Let $(R,s)\in\mathcal{Q}_{G,p}$. Then
$$\scalv{k}{F_{R,s}^G}_G=\left\{\begin{array}{cl}
\frac{\dsp\phi(|R|)}{\dsp|N_G(R,s)|}&\hbox{if ${<}sR{>}$ is cyclic}\\
0&\hbox{otherwise}
\end{array}\right.\mvirg$$
where $\phi$ is the Euler totient function.
\end{mth}
\pf Recall from \ref{idempotent} that
$$F_{R,{s}}^G=\frac{1}{|R||{s}||C_{\sur{N}_G(R)}({s})|}\sum_{\varphi,L}\tilde{\varphi}({s}^{-1})|L|\mu(L,{<}sR{>})\,\Ind_L^G\Res_L^{{<}sR{>}}k_\varphi\mvirg$$
where $\varphi$ runs through homomorphisms ${<}s{>}\to k^\times$ (i.e. equivalently homomorphisms ${<}sR{>}\to k^\times$), and $L$ through subgroups of ${<}sR{>}$ such that $LR={<}sR{>}$. Now 
$$\scalv{k}{\Ind_L^G\Res_L^{{<}sR{>}}k_\varphi}_G=\scalv{k}{k_\varphi}_L$$ 
is equal to zero unless the restriction of $\varphi$ to $L$ is trivial, i.e. if $L\leq\Ker\,\varphi$. Since $R\leq\Ker\,\varphi$, this implies that ${<}sR{>}\leq \Ker\,\varphi$, i.e. that $\varphi$ is trivial. Thus
$$\scalv{k}{F_{R,{s}}^G}_G=\frac{1}{|R||{s}||C_{\sur{N}_G(R)}({s})|}\sum_{LR={<}sR{>}}|L|\mu(L,{<}sR{>})\mpoint$$
Now a subgroup $L$ of $H={<}sR{>}$ is such that $LR=H$ if and only if it contains some conjugate of $s$ in $H$, i.e. if it is of the form $Q{<}^xs{>}$, for some $x\in R$ and some subgroup $Q$ of $R$ normalized by $^xs$. In this case, there are $|Q:C_Q({^xs})|$ conjugates of $s$ contained in $L$. Moreover, since $LR=H$ and $L\cap R=O_p(H)=Q$, the map $X\mapsto X\cap R$, from the poset $]L,H[$ of proper subgroups of $H$ strictly containing~$L$, to the poset $]Q,R[^{L}=]Q,R[^{^xs}$ of $L$-invariant proper subgroups of $R$ strictly containing~$Q$, is an isomorphism~: the inverse isomorphism is the map $Y\mapsto Y\cdot L$. Thus $\mu(L,{<}sR{>})$ is equal to the value $\mu\big((Q,R)^{^xs}\big)$ of the M\"obius function of the poset of subgroups normalized by $^xs$, and
\begin{eqnarray*}
\scalv{k}{F_{R,{s}}^G}_G\!&=&\!\!\frac{1}{|R||{s}||C_{\sur{N}_G(R)}({s})|}\!\sum_{x\in R/C_R(s)}\!\sumb{Q\leq R}{\rule{0cm}{2ex}Q^{^xs}=Q}\!|Q||s|\mu\big((Q,R)^{^xs}\big)/|Q:C_Q({^xs})|\\
&=&\frac{1}{|R||C_{\sur{N}_G(R)}({s})|}\sum_{x\in R/C_R(s)}\sumb{Q\leq R}{\rule{0cm}{2ex}Q^{^xs}=Q}|C_Q({^xs})|\mu\big((Q,R)^{^xs}\big)\mpoint
\end{eqnarray*}
Now $Q$ is normalized by $^xs$ if and only if $Q^x$ is normalized by $s$. Moreover $\mu\big((Q,R)^{^xs}\big)=\mu\big((Q^x,R)^{s}\big)$, and $|C_Q({^xs})|=|C_{Q^x}(s)|$. Thus
\begin{eqnarray*}
\scalv{k}{F_{R,{s}}^G}_G &=&\frac{1}{|R||C_{\sur{N}_G(R)}({s})|}\sum_{x\in R/C_R(s)}\sumb{Q\leq R}{\rule{0cm}{1.2ex}Q^{s}=Q}|C_Q({s})|\mu\big((Q,R)^{s}\big)\\
&=&\frac{1}{|C_R(s)||C_{\sur{N}_G(R)}({s})|}\sumb{Q\leq R}{\rule{0cm}{1.2ex}Q^{s}=Q}|C_Q({s})|\mu\big((Q,R)^{s}\big)\mpoint
\end{eqnarray*}
Now
\begin{eqnarray*}
\sumb{Q\leq R}{\rule{0cm}{1.2ex}Q^{s}=Q}|C_Q({s})|\mu\big((Q,R)^{s}\big)&=&\sum_{y\in C_R(s)}\sumb{y\in Q\leq R}{Q^s=Q}\mu\big((Q,R)^{s}\big)\\
&=&\sum_{y\in C_R(s)}\delta_{{<}y{>},R}\\
&=&|\{y\in C_R(s)\mid {<}y{>}=R\}|\mpoint
\end{eqnarray*}
This is non zero if and only if $R$ is cyclic and centralized by $s$, i.e. if the group ${<}sR{>}$ is cyclic. In this case $C_R(s)=R$, and
$$\scalv{k}{F_{R,{s}}^G}_G =\frac{1}{|R||C_{\sur{N}_G(R)}({s})|}\phi(|R|)=\frac{\phi(|R|)}{|N_G(R,s)|}\mvirg$$
since $|R||C_{\sur{N}_G(R)}({s})|=|N_G(R,s)|$, in view of Exact Sequence~\ref{suite exacte}. This completes the proof of the lemma.\findemo
\begin{mth}{Lemma}\label{sscal F} Let $(R,s)\in\mathcal{Q}_{G,p}$. Then
$$\sscalv{k}{F_{R,s}^G}_G=\frac{1}{|C_{\sur{N}_G(R)}(s)|}\sum_{x\in R/[{<}sR{>},R]}\frac{1}{|x|}\mpoint$$
\end{mth}
\pf By Definition~\ref{sscalv}
$$\sscalv{k}{F_{R,s}^G}_G=\sum_{Q\in[\mathcal{S}_p(G)]}\scalv{k}{F_{R,s}^G[Q]}_{\sur{N}_G(Q)}\mpoint$$
Hence, by Lemma~\ref{Brauer F} 
\begin{eqnarray*}
\sscalv{k}{F_{R,s}^G}_G&=&\sum_{Q\in\mathcal{S}_p(G)}\frac{|N_G(Q)|}{|G|}\sumb{x\in N_G(Q)\dom G/N_G(R,s)}{\rule{0ex}{1.5ex}Q\normal {^x{<}sR{>}}}\scalv{k}{F_{^xR/Q,^xs}^{\sur{N}_G(Q)}}_{\sur{N}_G(Q)}\\
&=&\sumc{Q\in\mathcal{S}_p(G)}{x\in G}{\rule{0ex}{1.5ex}Q^x\normal {<}sR{>}}\frac{|N_G(Q)|}{|G|}\frac{|N_G(Q^x)\cap N_G(R,s)|}{|N_G(Q)||N_G(R,s)|}\scalv{k}{F_{R/Q^x,s}^{\sur{N}_G(Q^x)}}_{\sur{N}_G(Q^x)}\\
&=&\sumb{Q\leq R}{\rule{0ex}{1.5ex}Q\normal {<}sR{>}}\frac{|N_G(Q)\cap N_G(R,s)|}{|N_G(R,s)|}\scalv{k}{F_{R/Q,s}^{\sur{N}_G(Q)}}_{\sur{N}_G(Q)}
\end{eqnarray*}
Now $\scalv{k}{F_{R/Q,s}^{\sur{N}_G(Q)}}_{\sur{N}_G(Q)}=0$ by Lemma~\ref{scal F}, unless $R/Q$ is cyclic and centralized by $s$, i.e. if $R/Q$ is cyclic and $Q$ contains the commutator subgroup $[{<}sR{>},R]$, in which case it is equal to $\frac{\dsp\phi(|R/Q|)}{\dsp|M|}$, where $M=N_{\sur{N}_G(Q)}(R/Q,s)$. Now the two exact sequences
$$\un\to R\to N_G(Q)\cap N_G(R,s)\to C_{N_G(Q,R)/R}(s)\to \un$$
$$\un\to R/Q\to N_{\sur{N}_G(Q)}(R/Q,s)\to C_{N_G(Q,R)/R}(s)\to \un$$
show that 
$$\frac{|N_G(Q)\cap N_G(R,s)|}{|N_{\sur{N}_G(Q)}(R/Q,s)|}=|Q|\mpoint$$
It follows that
\begin{eqnarray*}
\sscalv{k}{F_{R,s}^G}_G&=&\sumb{[{<}sR{>},R]\leq Q\leq R}{R/Q\;{\rm cyclic}}\frac{|Q|\phi(|R/Q|)}{|N_G(R,s)|}\\
&=&\frac{1}{|C_{\sur{N}_G(R)}(s)|}\sumb{[{<}sR{>},R]\leq Q\leq R}{R/Q\;{\rm cyclic}}\frac{\phi(|R/Q|)}{|R/Q|}\mpoint
\end{eqnarray*}
Now the summation is equal to the summation over all cyclic quotients $X=R/Q$ of the abelian group $A=R/[{<}sR{>},R]$ of $\frac{\dsp\phi(|X|)}{\dsp|X|}$, i.e. by a classical duality argument, to the summation of this quantity over the {\em cyclic subgroups} of  $A$, i.e. finally to $\dsp\sum_{x\in A}\limits\frac{1}{|x|}$. Thus
$$\sscalv{k}{F_{R,s}^G}_G=\frac{1}{|C_{\sur{N}_G(R)}(s)|}\sum_{x\in R/[{<}sR{>},R]}\frac{1}{|x|}\mvirg$$
as was to be shown.\findemo
\begin{mth}{Theorem} \label{theorem Mackey}Let $G$ be a finite group, let $p$ be a prime number, and $k$ be a field of characteristic~$p$, big enough to be a splitting field for all the groups $N_G(Q)/Q$, for $Q\in\mathcal{S}_p(G)$. Then the determinant of the Cartan matrix of the algebra $\mu_k(G,\un)$ is equal to
$$\det\mathsf{C}\big(\mu_k(G,\un)\big)=\prod_{R\in[\mathcal{S}_p(G)]}\prod_{s\in[\sur{N}_G(R)_{p'}]}\big(|C_{\sur{N}_G(R)}(s)|_p\sum_{x\in R/[{<}sR{>},R]}\frac{1}{|x|}\big)\mpoint$$
\end{mth}
\pf The Cartan matrix of $\mu_k(G,\un)$ is independent of the field $k$, as long as it is big enough. So one can assume e.g. that $k$ is algebraically closed, and choose a corresponding $p$-modular system $(K,\mathcal{O},k)$, where $K$ is also big enough. Then the formula for $\det\mathsf{C}\big(\mu_k(G,\un)\big)$ follows from Equation~\ref{detmu} and Lemma~\ref{sscal F}.\findemo
\begin{rem}{Examples} Recall that $k$ is a (big enough) field of characteristic $p$.\spn
$\bullet$ Let $G$ be a cyclic group of order $p^n$. Then
$$\det \mathsf{C}\big(\mu_k(G,\un)\big)=p^{\binom{n}{2}}\prod_{i=1}^n\big(p+i(p-1)\big)\mpoint$$
$\bullet$ Let $G$ be an elementary abelian group of order $p^2$. Then
$$\det \mathsf{C}\big(\mu_k(G,\un)\big)=p(2p-1)^{p+1}(p^2+p-1)\mpoint$$
\end{rem}
\section{The cohomological Mackey algebra}\label{coMackey}
\npar The cohomological Mackey algebra $co\mu_R(G)$ of a finite group $G$ over a commutative ring $R$ has also been introduced by Th\'evenaz and Webb (see~\cite{thevwebb} Section~16). It is an associative algebra with the property that the category of left $co\mu_R(G)$-modules is equivalent to the category $\comack_R(G)$ of {\em cohomological Mackey functors}. The algebra $co\mu_R(G)$ is defined as a quotient of the Mackey algebra $\mu_R(G)$.\par
In the case where $R$ is a field $k$ of characteristic $p$, the cohomological Mackey functors for $G$ over $k$ are projective relative to $p$-subgroups. In other words, the algebra $co\mu_k(G)$ is a quotient of $\mu_k(G,\un)$.
\npar Th\'evenaz and Webb have shown that the projective cohomological Mackey functors for $G$ over $R$ are exactly the fixed points functors $FP_W$, where $W$ is a direct summand of some permutation $RG$-module. Moreover, for any $RG$-modules $V$ and $W$,
$$\Hom_{\comack_R(G)}(FP_W,FP_V)\cong \Hom_{RG}(W,V)\mpoint$$
In particular, the indecomposable projective $co\mu_k(G)$-modules are the functors $FP_W$, where $W$ is an indecomposable $p$-permutation $kG$-module. Thus, the Cartan matrix of the algebra $co\mu_k(G)$ is the square matrix indexed by the set $[\mathcal{P}_{G,p}]$, defined by
$$c_{(P,E),(Q,F)}=\dim_k\Hom_{\comack_k(G)}(FP_{M_{P,E}},FP_{M_{Q,F}})=\scalv{M_{P,E}}{M_{Q,F}}_G\mpoint$$
\npar Thus by the same argument already used in Section~\ref{Mackey}, the Cartan matrix of the algebra $co\mu_k(G)$ is the matrix of the bilinear form $\scalv{~\,}{~}$ on $pp_k(G)$. This extends to a bilinear form on $K\otimes_\Z pp_k(G)$, and by~\ref{scalprod}, for $(R,s)$ and $(R',s')$ in $\mathcal{Q}_{G,p}$
$$\scalv{F_{R',s'}^G}{F_{R,s}^G}=\left\{\begin{array}{cl}
\scalv{k}{F_{R,s}^G}_G&\hbox{if}\;(R',s')=_G(R,s^{-1})\\
0&\hbox{otherwise}\end{array}\right.$$
By Lemma~\ref{scal F}, this is non zero if and only if the group ${<}sR{>}$ is cyclic.\par
The argument of paragraph~\ref{a suivre} shows that the Cartan matrix of $co\mu_k(G)$ can be expressed as
$$\mathsf{C}\big(co\mu_k(G)\big)=\mathsf{T}'\cdot \mathsf{S'}\cdot{^t\mathsf{T}}\mvirg$$
where $\mathsf{S}'$ is the diagonal matrix indexed by $[\mathcal{Q}_{G,p}]$ with $(R,s)$-diagonal entry equal to $\scalv{F_{R,s}^G}{k}_G$. Since $\mathsf{T}$ and $\mathsf{T}'$ are invertible, it follows that the rank of $\mathsf{C}\big(co\mu_k(G)\big)$ is equal to the number of elements $(R,s)$ of $[\mathcal{Q}_{G,p}]$ such that ${<}sR{>}$ is cyclic.
\begin{mth}{Theorem} \label{theorem coMackey}Let $G$ be a finite group, let $p$ be a prime number, and $k$ be a field of characteristic~$p$, big enough to be a splitting field for all the groups $N_G(Q)/Q$, for $Q\in\mathcal{S}_p(G)$. Then~:
\begin{enumerate}
\item The rank of the Cartan matrix of the algebra $co\mu_k(G)$ is equal to the number of $G$-conjugacy classes of pairs $(R,s)$, where $R\in\mathcal{S}_p(G)$ and $s\in\sur{N}_G(R)_{p'}$, such that ${<}sR{>}$ is cyclic. In other words
\begin{eqnarray*}
{\rm rk}\,\mathsf{C}\big(co\mu_k(G)\big)&=&\sum_{R\in [\mathcal{C}_p(G)]}|N_G(R)\dom C_G(R)_{p'}|\\
&=&\sum_{s\in[G_{p'}]}c_p\big(C_G(s)\big)\mvirg
\end{eqnarray*}
where $[\mathcal{C}_p(G)]$ denotes a set of representatives of conjugacy classes of cyclic $p$-subgroups of a group $G$, and $c_p(G)$ is the number of such conjugacy classes. Moreover, in the first summation, $|N_G(R)\dom C_G(R)_{p'}|$ denotes the number of $N_G(R)$-conjugacy classes of $p$-regular elements of $C_G(R)$.
\item The Cartan matrix of the algebra $co\mu_k(G)$ is non singular if and only if the group $G$ is $p$-nilpotent with cyclic Sylow $p$-subgroups. In this case, if $G=N\rtimes P$, where $N$ is a $p'$-group and $P$ is a cyclic group of order~$p^n$
\begin{eqnarray*}
{\rm det}\,\mathsf{C}\big(co\mu_k(G)\big)&=&\prod_{R\leq P}\left(\frac{\phi(|R|)}{|R|}\right)^{l_p\big(\sur{N}_G(R)\big)}{\rm det}\,\mathsf{C}\big(k\sur{N}_G(R)\big)\\
&=&\prod_{R\leq P}\prod_{x\in \big[P\dom [C_N(R)]\big]}\left(\frac{\phi(|R|)}{|R|}|C_{C_N(R)\rtimes P/R}(x)|_p\right)\mvirg
\end{eqnarray*}
where $l_p\big(\sur{N}_G(R)\big)$ is the number of $p$-regular classes of $\sur{N}_G(R)$, where $\mathsf{C}\big(k\sur{N}_G(R)\big)$ is the Cartan matrix of the group algebra $k\sur{N}_G(R)$, and $\big[P\dom [C_N(R)]\big]$ is a set of representatives of $C_N(R)\rtimes P$-conjugacy classes in $C_N(R)$.
\end{enumerate}
\end{mth}
\pf The first sentence of Assertion~1 follows from the above arguments, and the second one follows by counting the number of pairs $(R,s)$ of $[\mathcal{Q}_{G,p}]$ such that ${<}sR{>}$ is cyclic, in two different ways. For Assertion~2, observe that the Cartan matrix $\mathsf{C}\big(co\mu_k(G)\big)$ is non singular if and only if for any $(R,s)\in\mathcal{Q}_{G,p}$, the group ${<}sR{>}$ is cyclic. This amounts to saying that the Sylow $p$-subgroups of $G$ are cyclic, and that whenever $s\in G_{p'}$ normalizes a $p$-subgroup $R$, it centralizes it. In other words, the group $N_G(R)/C_G(R)$ is a $p$-group, for any $R\in\mathcal{S}_p(G)$. This is equivalent to saying that $G$ is $p$-nilpotent, by the theorem of Frobenius (\cite{gorenstein} Theorem 4.5).\par
In the case $G=N\rtimes P$, where $P$ is cyclic of order $p^n$, then as in~\ref{detmu}
\begin{eqnarray*}
\det \mathsf{C}\big(co\mu_k(G)\big)\!\!&\!=\!&\!\det\mathsf{T}'\,\det\mathsf{S'}\,\det{^t\mathsf{T}}\\
\!\!&\!=\!&\!\prod_{R\in[\mathcal{S}_p(G)]}\prod_{{s}\in[\sur{N}_G(R)_{p'}]}\big(|C_{\sur{N}_G(R)}({s})|_p|C_{\sur{N}_G(R)}({s})|\scalv{k}{F_{R,{s}}^G}_G\big)\\
\!\!&\!=\!&\!\prod_{R\in[\mathcal{S}_p(G)]}\!\left(\!\det \mathsf{C}\big(k\sur{N}_G(R)\big)\!\!\!\!\prod_{{s}\in[\sur{N}_G(R)_{p'}]}\!\!\!\!\big(|C_{\sur{N}_G(R)}({s})|\scalv{k}{F_{R,{s}}^G}_G\big)\!\right).\\
\end{eqnarray*}
The set $[\mathcal{S}_p(G)]$ can be chosen to be the set of subgroups of $P$. Moreover by Lemma~\ref{scal F}, for $R\leq P$ and $s\in\sur{N}_G(R)_{p'}$
$$ \scalv{k}{F_{R,{s}}^G}_G=\frac{\phi(|R|)}{|N_G(R,s)|}=\frac{\phi(|R|)}{|R||C_{\sur{N}_G(R)}(s)|}\mpoint
$$
It follows that
$$\det\,\mathsf{C}\big(co\mu_k(G)\big)=\prod_{R\leq P}\left(\frac{\phi(|R|)}{|R|}\right)^{l_p\big(\sur{N}_G(R)\big)}{\rm det}\,\mathsf{C}\big(k\sur{N}_G(R)\big)\mpoint$$
Finally the $p'$-elements of the group $\sur{N}_G(R)\cong C_N(R)\rtimes (P/R)$ are the elements of $C_N(R)$. The last formula of the theorem follows.\findemo 
\begin{rem}{Remark} Thus when it is non zero, the determinant of the Cartan matrix of $co\mu_k(G)$ is equal to $(p-1)^np^m$, for suitable non negative integers $n$ and~$m$.
\end{rem}
\section{Blocks of cohomological Mackey functors}\label{coMackey blocks}
\npar It was shown by Th\'evenaz and Webb (see \cite{thevwebb} Theorem~17.1 and its proof, see also~\cite{mublocks})) that the blocks of the algebra $\mu_k(G,\un)$ and the blocks of the algebra $co\mu_k(G)$ are in one to one correspondence with the blocks of the group algebra $kG$. If $b$ is a block of $kG$, denote by $co\mu_k(b)$ the corresponding block of $co\mu_k(G)$. When $R$ is a $p$-subgroup of $G$, let ${\rm Br_R}:(kG)^R\to kC_G(R)$ denote the Brauer morphism. If moreover $R\normal G$, denote by $u\mapsto \sur{u}$ the projection map $kG\to k(G/R)$. When $A$ is a $k$-algebra, denote by $\Irr_k(A)$ the set of isomorphism classes of simple $A$-modules.\par
This section is devoted to the proof of the following block version of Theorem~\ref{theorem coMackey} for such blocks of cohomological Mackey functors~:
\begin{mth}{Theorem} \label{theorem coMackey blocks}Let $G$ be a finite group, let $p$ be a prime number, and $k$ be an algebraically closed field of characteristic~$p$. Let moreover $b$ be a block of $kG$. Then~: 
\begin{enumerate}
\item The rank of the Cartan matrix of the algebra $co\mu_k(b)$ is equal to
\begin{eqnarray*}
{\rm rk}\;\mathsf{C}\big(co\mu_k(b)\big)&=&\sum_{R\in[\mathcal{C}_p(G)]}|N_G(R)\dom{\rm Irr}_k\big(kC_G(R){\rm Br}_R(b)\big)|\\
&=&\sum_{(R,c)\in [\mathcal{C}_p(b)]}|N_G(R,c)\dom{\rm Irr}_k\big(kC_G(R)c\big)|\mvirg
\end{eqnarray*}
where $[\mathcal{C}_p(b)]$ is a set of representatives of $G$-conjugacy classes of $b$-Brauer pairs $(R,c)$ for which $R$ is cyclic.
\item The Cartan matrix of the block $co\mu_k(b)$ is non singular if and only if $b$~is a nilpotent block with cyclic defect groups.
\end{enumerate}
\end{mth}
\npar The Cartan matrix $\mathsf{C}$ of $co\mu_k(b)$ is non singular if and only if the rows of the decomposition matrix $\mathsf{D}$ are linearly independent (since $\mathsf{C}=\mathsf{D}\cdot{^t\mathsf{D}}$)~: indeed, if a vector $u$ is such that $\mathsf{C}u=0$, then ${^tu}\mathsf{C}u={^t({^t\mathsf{D}}u)}\cdot({^t\mathsf{D}}u)=0$, thus ${^t\mathsf{D}}u=0$. Conversely, if ${^t\mathsf{D}}u=0$, then obviously $\mathsf{C}u=0$.
\npar The decomposition matrix $\mathsf{D}$ of $co\mu_k(b)$ has rows indexed by the (isomorphism classes of) indecomposable $p$-permutation $kG$-modules in the block~$b$, and columns indexed by the (isomorphism classes of) simple $KG$-modules in the block $b$. The coefficient $\mathsf{D}_{W,\chi}$ corresponding to the indecomposable $p$-permutation $kGb$-module $W$ and the simple $KG$-module $\chi$ is equal to the multiplicity of $\chi$ in $K\otimes_\mathcal{O}\tilde{W}$, where $\tilde{W}$ is an $\mathcal{O}G$-module lifting $W$ to $\mathcal{O}$ (i.e. such that $k\otimes_\mathcal{O}\tilde{W}\cong W$). Such an $\mathcal{O}G$-module is unique up to isomorphism, since $W$ is a $p$-permutation module. The character of the module $K\otimes_\mathcal{O}\tilde{W}$ will be called {\em the (ordinary) character} of $W$. \par
It follows that the rank of the Cartan matrix of $co\mu_k(b)$ is equal to the dimension of the subspace of $\Q\otimes_\Z R_K(b)$ generated by characters of $p$-permutation modules in the block $b$.\par
In particular, the Cartan matrix of $co\mu_k(b)$ is non singular if and only if the ordinary characters of the indecomposable $p$-permutation modules in $b$ are linearly independent.
\npar \label{value}Recall (see~\ref{Brauer trace}, and \cite{broue-galway} Proposition 3.3) that the value of the ordinary character of a $p$-permutation module $W$ on an element $s$ of $G$ is equal to the Brauer trace ${\rm BrTr}(s_{p'}\mid W[{<}s_p{>}])$, where $s_p$ and $s_{p'}$ are the $p$-part and $p'$-part of $s$, respectively. 
\begin{rem}{Notation} Let $\mathcal{Z}_p(G)$ denote the set of pairs $(P,E)$ consisting of a cyclic $p$-subgroup $P$ of $G$ and an indecomposable projective $k\sur{C}_G(P)$-module~$E$, where $\sur{C}_G(P)=C_G(P)/P$, and let $[\mathcal{Z}_p(G)]$ be a set of representatives of $G$-orbits on $\mathcal{Z}_p(G)$.\par
Let $\mathcal{Z}_p(b)$ denote the subset of $\mathcal{Z}_p(G)$ consisting of pairs $(P,E)$ such that 
$${\rm Br}_P(b)\Inf_{\sur{C}_G(P)}^{C_G(P)}E=\Inf_{\sur{C}_G(P)}^{C_G(P)}E\mpoint$$
Set moreover $[\mathcal{Z}_p(b)]=\mathcal{Z}_p(b)\cap [\mathcal{Z}_p(G)]$.
\end{rem}
\begin{mth}{Lemma} \label{characters}The characters of the modules $\Ind_{C_G(P)}^G\Inf_{\sur{C}_G(P)}^{C_G(P)}E$, for $(P,E)\in [\mathcal{Z}_p(G)]$, form a basis of the subspace of $\Q\otimes_\Z R_K(G)$ generated by the characters of the $p$-permutation modules. 
\end{mth}
\pf Since projective modules are $p$-permutation modules, and since induction and inflation preserves this class of modules, the module $L_{P,E}=\Ind_{C_G(P)}^G\Inf_{\sur{C}_G(P)}^{C_G(P)}E$, for $(P,E)\in\mathcal{Z}_p(G)$, is a $p$-permutation module. Up to isomorphism, this module depends only on the $G$-orbit of $(P,E)$. Moreover the number of $G$-orbits on the set $\mathcal{Z}_p(G)$ is equal to 
$$|G\dom\mathcal{Z}_p(G)|=\sum_{P\in[\mathcal{C}_p(G)]}|N_G(P)\dom {\rm Irr}_k\big(kC_G(P)\big)|\mpoint$$
Indeed, the number of isomorphism classes of indecomposable projective $k\sur{C}_G(P)$-modules is equal to the number of isomorphism classes of simple $k\sur{C}_G(P)$-modules, i.e. to the number of isomorphism classes of simple $kC_G(P)$-modules, since $P$ is a normal $p$-subgroup of $C_G(P)$.\par
By Theorem~\ref{theorem coMackey}, it follows that the cardinality of the set $[\mathcal{Z}_p(G)]$ is precisely equal to the rank of the Cartan matrix of the algebra $co\mu_k(G)$, i.e. to the dimension of the subspace of $\Q\otimes_\Z R_K(G)$ generated by the characters of the $p$-permutation modules. Thus, to prove Lemma~\ref{characters}, it is enough to prove that the characters of the modules $L_{P,E}$, for $(P,E)\in[\mathcal{Z}_p(G)]$, are linearly independent.\par
The character $\chi_{P,E}$ of the module $L_{P,E}$ is equal to $\Ind_{C_G(P)}^G\Inf_{\sur{C}_G(P)}^{C_G(P)}\Phi_E$, where $\Phi_E$ is the character of the module $E$. Suppose that some non trivial linear combination of these characters is equal to 0, i.e. that there are integers $n_{P,E}\in \Z$, for $(P,E)\in[\mathcal{Z}_p(G)]$, not all equal to 0, such that
\begin{equation}\label{combination}
\sum_{(P,E)\in[\mathcal{Z}_p(G)]}n_{P,E}\chi_{P,E}=0\mpoint
\end{equation}
Let $Q$ be maximal such that there exists $(Q,F)\in [\mathcal{Z}_p(G)]$ with $n_{Q,F}\neq 0$, and let $s$ be a generator of the cyclic group $Q$. By~\ref{value}, for any $t\in C_G(s)_{p'}$, the value of the linear combination~\ref{combination} at the element $st$ is equal to 
$$0=\sum_{(P,E)\in[\mathcal{Z}_p(G)]}n_{P,E}{\rm BrTr}(t\mid L_{P,E}[Q])\mpoint$$
But $L_{P,E}[Q]=0$, unless some conjugate of $Q$ is contained in $P$. By maximality of~$Q$, it follows that
$$0=\sum_{E}n_{Q,E}{\rm BrTr}(t\mid L_{Q,E}[Q])\mvirg$$
where $E$ runs through a set $[\mathcal{P}]$ of indecomposable projective $k\sur{C}_G(Q)$-modules, up to isomorphism and conjugation by $N_G(Q)$. Since moreover $$L_{Q,E}[Q]\cong \Ind_{C_G(Q)}^{N_G(Q)}\Inf_{\sur{C}_G(Q)}^{C_G(Q)}E\mvirg$$
it follows that for any $t\in C_G(Q)_{p'}$
$$\sum_{E\in[\mathcal{P}]}n_{Q,E}\big(\Ind_{C_G(Q)}^{N_G(Q)}\Inf_{\sur{C}_G(Q)}^{C_G(Q)}\Phi_E\big)(t)=0\mpoint$$
This is also equal to
\begin{eqnarray*}
\sum_{E\in[\mathcal{P}]}n_{Q,E}\big(\Res_{C_G(Q)}^{N_G(Q)}\Ind_{C_G(Q)}^{N_G(Q)}\Inf_{\sur{C}_G(Q)}^{C_G(Q)}\Phi_E\big)(t)&=&\sumb{E\in[\mathcal{P}]}{x\in N_G(Q)/C_G(Q)}n_{Q,E}\Phi_E({\sur{^xt}})\\
&=&\sumb{E\in[\mathcal{P}]}{x\in N_G(Q)/C_G(Q)}n_{Q,E}\Phi_{E^x}({\sur{t}})\mvirg
\end{eqnarray*}
where $\sur{t}$ is the image of $t$ in $\sur{C}_G(Q)$, and $E^x$ is the image of $E$ by conjugation by $x\in N_G(Q)$.\par
Since the map $t\mapsto \sur{t}$ is a surjection from $C_G(Q)_{p'}$ to $\sur{C}_G(Q)_{p'}$, it follows that the modular character $\sumb{E\in[\mathcal{P}]}{x\in N_G(Q)/C_G(Q)}\limits n_{Q,E}\Phi_{E^x}$ of $\sur{C}_G(Q)$ is equal to zero. But the set of modules $E^x$, for $E\in[\mathcal{P}]$ and $x\in N_G(Q)/C_G(Q)$, is exactly the set of projective indecomposable $k\sur{C}_G(Q)$-modules, up to isomorphism. Thus
$$\sum_{F\in\mathcal{Q}}n_{Q,F}m_F\Phi_F=0\mvirg$$
where $\mathcal{Q}$ is a set of representatives of isomorphism classes of indecomposable projective $k\sur{C}_G(Q)$-modules, where $n_{Q,F}$ is defined as $n_{Q,E}$ if $(Q,E)\in [\mathcal{P}]$ and if there exists $x\in N_G(Q)$ such that $E^x\cong F$, and where $m_F$ is the number of elements $x\in N_G(Q)/C_G(Q)$ such that $F^x\cong F$. \par
Now the characters $\Phi_F$, for $F\in\mathcal{Q}$, are linearly independent, and it follows that $n_{Q,F}=0$ for any $F$. This contradicts the definition of $Q$, and completes the proof of Lemma~\ref{characters}.\findemo
\begin{rem}{Notation} Let $P\in\mathcal{S}_p(G)$, and $E$ be any projective $k\sur{N}_G(P)$-module. Then $E$ splits as a direct sum $\mathop{\oplus}_{i\in I}\limits E_i$ of indecomposable $k\sur{N}_G(P)$-modules $E_i$. In this situation, set $M_{P,E}=\mathop{\oplus}_{i\in I}\limits M_{P,E_i}$.
\end{rem}
\begin{mth}{Lemma} \label{characters-blocks}The characters of the modules $M_{P,\Ind_{\sur{C}_G(P)}^{\sur{N}_G(P)}E}$, for $(P,E)\in[\mathcal{Z}_p(b)]$, form a basis of the subspace of $\Q\otimes_\Z R_K(b)$ generated by the characters of the $p$-permutation modules in the block $b$.
\end{mth}
\pf Let $(P,E)\in \mathcal{Z}_p(b)$. The module $M=\Ind_{\sur{C}_G(P)}^{\sur{N}_G(P)}E$ splits as a direct sum $\mathop{\oplus}_{i\in I}\limits E_i$ of indecomposable $k\sur{N}_G(P)$-modules $E_i$. Since ${\rm Br}_P(b)$ is $N_G(P)$-invariant and belongs to $kC_G(P)$, it follows that $\sur{{\rm Br}_P(b)}$ acts as the identity on $M$, hence on every direct summand~$E_i$. Recall that the indecomposable module $M_{P,E_i}$ belongs to the block $b$ if and only if $\sur{{\rm Br}_P(b)}E_i=E_i$ (see \cite{benson1}~Corollary 6.3.2). Hence all the indecomposable $kG$-modules $M_{P,E_i}$ are in the block~$b$, and their direct sum $M_{P,\Ind_{\sur{C}_G(P)}^{\sur{N}_G(P)}E}$ is also in $b$.\par
To prove Lemma~\ref{characters-blocks}, it suffices to observe that the sets $\mathcal{Z}_p(b)$, when $b$ runs through the blocks of $kG$, form a partition of $\mathcal{Z}_p(G)$, and to prove that the characters of the modules $M_{P,\Ind_{\sur{C}_G(P)}^{\sur{N}_G(P)}E}$, for $(P,E)\in[\mathcal{Z}_p(G)]$, form a basis of the subspace of $\Q\otimes_\Z R_K(G)$ generated by the characters of the $p$-permutation modules.\par
For this, observe that
$$\Ind_{C_G(P)}^G\Inf_{\sur{C}_G(P)}^{C_G(P)}E\cong \Ind_{N_G(P)}^G\Inf_{\sur{N}_G(P)}^{N_G(P)}\Ind_{\sur{C}_G(P)}^{\sur{N}_G(P)}E$$
is a direct sum of $M_{P,\Ind_{\sur{C}_G(P)}^{\sur{N}_G(P)}E}$ and of indecomposable modules with vertex strictly contained in $P$ up to conjugation. This yields a triangular transition matrix, with non zero diagonal coefficients, and such a matrix changes a basis to another basis.\findemo
\npar  The first equality in Assertion 1 of Theorem~\ref{theorem coMackey blocks} follows trivially from Lemma~\ref{characters-blocks}. The second one follows from the fact that, for any $b$-Brauer pair $(R,c)$
$$|N_G(R)\dom {\rm Irr}_k\big(kC_G(R){\rm Tr}_{N_G(R,c)}^{N_G(R)}c\big)|=|N_G(R,c)\dom {\rm Irr}_k\big(kC_G(R)c\big)|\mpoint$$
This is because the algebra $kC_G(R){\rm Tr}_{N_G(R,c)}^{N_G(R)}c$ is isomorphic to the direct sum of the block algebras $kC_G(R)\,{^xc}$, for $x\in [N_G(R)/N_G(R,c)]$, which are transitively permuted by $N_G(R)$. Moreover the stabilizer in $N_G(R)$ of $kC_G(R)c$ is equal to $N_G(R,c)$. 
\npar For any $p$-subgroup $R$ of $G$, induction from $C_G(R)$ to $N_G(R)$ induces an inequality
\begin{equation}\label{inequality}
|N_G(R)\dom {\rm Irr}_k\big(kC_G(R){\rm Br}_R(b)\big)|\leq |{\rm Irr}_k\big(kN_G(R){\rm Br}_R(b)\big)|\mpoint
\end{equation}
Indeed, the left hand side is equal to the rank of the group
$$P_k\big(kC_G(R){\rm Br}_R(b)\big)_{N_G(R)}$$
of $N_G(R)$-coinvariants on the group of projective $kC_G(R){\rm Br}_R(b)$-modules. It is easy to see that $\Ind_{C_G(R)}^{N_G(R)}$ induces an injective map from this group to the group $P_k\big(kN_G(R){\rm Br}_R(b)\big)$ of projective $kN_G(R){\rm Br}_R(b)$-modules, and the rank of this group is equal to the right hand side of~\ref{inequality}.\par
The Cartan matrix of $co\mu_k(b)$ is non singular if and only if
$$\sum_{R\in[\mathcal{C}_p(G)]}|N_G(R)\dom {\rm Irr}_k\big(kC_G(R){\rm Br}_R(b)\big)|=\sum_{R\in[\mathcal{S}_p(G)]}|{\rm Irr}_k\big(kN_G(R){\rm Br}_R(b)\big)|\mpoint$$
Indeed, the left hand side is the rank of the Cartan matrix, and the right hand side the size of this matrix, i.e. the number of indecomposable $p$-permutation modules in the block $b$, up to isomorphism.\par
By inequality~\ref{inequality}, this is in turn equivalent to the following equality, for any $R\in\mathcal{S}_p(G)$~:
$$|{\rm Irr}_k\big(kN_G(R){\rm Br}_R(b)\big)|=\left\{\begin{array}{cl}|N_G(R)\dom {\rm Irr}_k\big(kC_G(R){\rm Br}_R(b)\big)|&\hbox{if $R$ is cyclic}\\0&\hbox{otherwise}\mpoint\end{array}\right.$$
\npar Hence if the Cartan matrix of $co\mu_k(b)$ is non singular, and $b$ has defect $D$, then in particular ${\rm Br}_D(b)\neq 0$, and $|{\rm Irr}_k\big(kN_G(D){\rm Br}_D(b)\big)|\neq 0$, so $D$ is cyclic. Let $(D,c)$ be a maximal $b$-Brauer pair. Then ${\rm Br}_D(b)={\rm Tr}_{N_G(D,c)}^{N_G(D)}c$, and the algebra $kN_G(D){\rm Br}_D(b)$ is isomorphic to a matrix algebra over $kN_G(D,c)c$. In particular
\begin{equation}\label{defect}
|{\rm Irr}_k\big(kN_G(D){\rm Br}_D(b)\big)|=|{\rm Irr}_k\big(kN_G(D,c)c\big)|\mpoint
\end{equation}
Moreover ${\rm Br}_D(b)$ splits as a sum of blocks of $C_G(D)$, which are the distinct $N_G(D)$-conjugates of $c$. Hence, the orbits of $N_G(D)$ on ${\rm Irr}_k\big(kC_G(D){\rm Br}_D(b)\big)$ are in one to one correspondence with the orbits of $N_G(D,c)$ on ${\rm Irr}_k\big(kC_G(D)c\big)$. It follows that
$$|N_G(D,c)\dom{\rm Irr}_k\big(kC_G(D)c\big)|=|{\rm Irr}_k\big(kN_G(D,c)c\big)|\mpoint$$
Now the group $N_G(D,c)/C_G(D)$ is the inertial quotient of $b$. It is a cyclic group of order $e$ dividing $p-1$. The block $c$ is a nilpotent block of $C_G(D)$, so in particular there is a unique simple $kC_G(D)c$-module $S$, which is invariant by $N_G(D,c)$. This simple module can be extended to a simple $kN_G(D,c)c$ module in $e$ different ways, by the following argument (see \cite{benson1}, proof of Proposition~6.5.4)~: let $g$ be a generator of the group $N_G(Q)/C_G(Q)$, and $\theta: S\to S^g$ an isomorphism of $kC_G(Q)$-modules. Then the map $g^e(s)\mapsto \theta^e(s)$ is a $kC_G(Q)$ automorphism of $S$, because $g^e\in C_G(Q)$. As $k$ is algebraically closed, there is a scalar $\mu\in k$ such that $g^e(s)=\mu\theta^e(s)$, for any $s\in S$. Moreover there are $e$ distinct elements $\lambda\in k$ such that $\lambda^e=\mu$. For each such~$\lambda$, one can let $g$ act on $S$ by $\lambda\theta$, and this gives $e$ mutually non isomorphic extensions of $S$ to a simple $kN_G(Q,c)$-module, which are all in the block $c$ since $c\in kC_G(Q)$. \par
These modules are not isomorphic to each other (since their restrictions to ${<}g{>}$ are not). Hence $e=1$ by~\ref{defect}, since $|{\rm Irr}_k\big(kC_G(D)c\big)|=1$, and this is equivalent to saying that $b$ is nilpotent, since $b$ has cyclic defect.
\npar Conversely, if $b$ is a nilpotent block with a cyclic defect group $D$, then for each $R\in\mathcal{S}_p(G)$, either ${\rm Br}_R(b)=0$ if $R$ is not contained in $D$ up to $G$ conjugation, or ${\rm Br}_R(b)$ is a sum $\sum_{i\in I}\limits b_i$ of blocks $b_i$ of $N_G(R)$, if $R$ is contained in $D$ up to conjugation, and in that case $R$ is cyclic since $D$ is. Each of the blocks $b_i$ is equal to ${\rm Tr}_{N_G(R,c_i)}^{N_G(R)}c_i$, where $(R,c_i)$ is some $b$-Brauer pair. In particular
\begin{eqnarray*}
|{\rm Irr}_k\big(kN_G(R){\rm Br}_R(b)\big)|&=&\sum_{i\in I} |{\rm Irr}_k\big(kN_G(R){\rm Br}_R(b_i)\big)|\\
&=&\sum_{i\in I}|{\rm Irr}_k\big(kN_G(R,c_i)c_i\big)|\mpoint
\end{eqnarray*}
Since the block $b$ is nilpotent, all the groups $N_G(R,c_i)/C_G(R)$ are $p$-groups, so $|{\rm Irr}_k\big(kN_G(R,c_i)c_i\big)|=|N_G(R,c_i)\dom {\rm Irr}_k\big(kC_G(R)c_i\big)|$, and this is equal to~1 since $c_i$ is a nilpotent block of $C_G(R)$ by Theorem~1.2 of~\cite{broue-puig_frobenius}. It follows that 
$$|{\rm Irr}_k\big(kN_G(R){\rm Br}_R(b)\big)|=|I|=|N_G(R)\dom {\rm Irr}_k\big(kC_G(R){\rm Br}_R(b)\big)|\mpoint$$
This implies that the Cartan matrix of $co\mu_k(b)$ is non singular, and completes the proof of Theorem~\ref{theorem coMackey blocks}.\findemo

\begin{flushleft}
Serge Bouc\\
LAMFA - CNRS UMR 6140\\
Universit\'e de Picardie Jules Verne, 33, rue St Leu\\
80039 - Amiens - FRANCE.\\
{\tt email : serge.bouc@u-picardie.fr}\\
{\tt web : http://www.mathinfo.u-picardie.fr/bouc/}
\end{flushleft}
\end{document}